\documentclass[a4paper,11pt]{article}
\usepackage{amssymb,amsmath,amsthm}
\usepackage{bm}
\usepackage{color}
\usepackage[normalem]{ulem}

\newtheorem{theorem}{Theorem}[section]
\newtheorem{proposition}[theorem]{Proposition}
\newtheorem{remark}[theorem]{Remark}
\newtheorem{lemma}[theorem]{Lemma}
\newtheorem{corollary}[theorem]{Corollary}
\newtheorem{definition}[theorem]{Definition}

\begin{document}

\title{Solvability of the asymmetric Bingham fluid equations}
\author{}
\maketitle

\centerline{\scshape Anderson Luis Albuquerque de Araujo} \medskip 
{\footnotesize \ 
\centerline{Departamento de Matemática,} 
\centerline{Universidade Federal
de Vi\c{c}osa. Vi\c{c}osa, MG, Brasil} 
\centerline{ Email:
anderson.araujo@ufv.br} }

\centerline{\scshape Nikolai V. Chemetov} \medskip {\footnotesize 
\centerline{Universidade de Lisboa, Edificio C6, 1 Piso, Campo Grande,} %
\centerline{1749-016 Lisboa, Portugal} 
\centerline{ Email:
nvchemetov@fc.ul.pt, \ \ nvchemetov@gmail.com}}

\medskip

\centerline{\scshape Marcelo M. Santos } \medskip {\footnotesize 
\centerline{Departamento de Matem\'atica,} 
\centerline{IMECC-Instituto de Matem\'atica, Estat\'\i stica
e Computa\c c\~ao Cient\'\i fica,} 
\centerline{UNICAMP-Universidade Estadual
de Campinas. Campinas, SP, Brazil} 
\centerline{ Email:
msantos@ime.unicamp.br  }}

\begin{abstract}
In this work, we investigate the asymmetric Bingham fluid equations. The asymmetric fluid of Bingham includes symmetric and antisymmetric stresses with such stresses appearing as an elastic response to the micro-rotational deformations of grains in a complex fluid. We show the global-in-time solvability of a weak solution for three dimensional boundary value problem with Navier boundary conditions of the asymmetric Bingham fluid equations.
\end{abstract}

\textit{Keywords}: asymmetric fluids, Bingham plastic fluids, yield stress,
global existence, weak solution.

MSC (2010):  35Q35, 74A20, 76S05

\bigskip

\section{Introduction}

\label{Introduction}

\bigskip

The Newtonian fluid obeys the constitutive relations that the deviation
stress tensor is a linear function of the stress rate-of-strain tensor and
if the fluid is isotropic then the stress and rate-of strain tensors are
symmetric. \ Many of fluids can not be described by the Newtonian
constitutive relations, such as slurries, animal blood, mud (mixtures of
water, clay), viscous polymers, polymeric suspensions. These fluids does not
commence to flow till the applied stress attains a certain optimum
magnitude, called the \textit{yield stress }$\tau _{\ast }$, after of that
they behave as a Newtonian fluid. \ An example is toothpaste, which will not
be extruded until a certain pressure is applied to the tube. Then it is
pushed out as a solid plug. The physical reason for such behaviour is that
the fluid contains particles, such as toothpaste, paints, clay, or large
molecules, such as polymers, animal blood, which have an interaction,
creating a weak rigid structure. Therefore a certain stress is required to
break this weak rigid structure. As soon as the structure has been broken,
the particles move with the fluid under viscous forces. The particles will
associate again if the stress is removed. Such behaviour was firstly
presented in an experimental study by Bingham \cite{B}, where he proposed
its mathematical model. Later on these type of fluids have been called as
Bingham plastic fluids. The Bingham plastic fluid behaves as a rigid body at
low stresses but flows as a viscous fluid at high stress. Nowadays it is
used as a mathematical model of mud flows in drilling engineering, heavy
oil, lava (being a mix with melting snow, stones), in the handling of
slurries, waxy crude oils. Recent examples concerns the propane flow within
the hydro-fracture \cite{ShelNeverov2014}. Significant efforts of the study
of the Bingham plastic fluids have been done by Oldroyd \cite{ol}, Mossolov,
Miasnikov \cite{MM}, Glowinski, Wachs \cite{GW}, Papanastasiou \cite%
{Papanastasiou} and many others scientists.

The Newtonian flow in the Bingham fluid (after the load greater than the 
\textit{yield stress }$\tau _{\ast }$) has a drastic limitation, since it
does not account the behaviour of the fluid, that contains the particles.
Most of above mentioned fluid systems contain rigid, randomly oriented
particles, irregularly shaped particles (drops in emulsions), branched and
entangled molecules in case of polymeric systems, or loosely formed clusters
of particles in suspensions, etc. The particles may shrink and expand or
change their shape, they may \textit{rotate}, independently of the rotation
of the fluid. To describe accurately the behaviour of such fluids a
so-called \textit{asymmetric continuum theory} \cite{T} (or \textit{\
micropolar theory} \cite{C}, \cite{E}) has been developed that ignores the
deformation of the particles but takes into account geometry, intrinsic
motion of material particles. \ This theory is a significant and a simple
generalization of the classical Navier-Stokes model, that describe the
Newtonian fluids. Only one new vector field, called as the \textit{angular
velocity} field, of rotation of particles is introduced. As a consequence,
only one equation is added, that represents the conservation of the angular
momentum. \ The asymmetric/micropolar fluids\textit{\ }belong to the class
of fluids with non-symmetric stress tensor. This class of fluids is more
general \ than the classical Newtonian fluids.

Shelukhin, R\r{u}\v{z}i\v{c}ka \cite{ShelRuz2013} suggested a mathematical
model that describe the behaviour of an asymmetric/micropolar viscous fluid
of Bingham. \ In the paper \cite{ShelChem} the global solvability for the
mathematical model of \cite{ShelRuz2013} has been demonstrated in a special
case of one dimensional flow. Later on in \cite{smr} the authors have
obtained the solvability of the stationary solution for this model already
for three dimensional case, but in a particular case when the stress tensor
does not have non-symmetric part.

The main objective of the current work is to correct the model suggested in 
\cite{ShelRuz2013} and to show the well-posedness of a modified model.

\bigskip

The paper is organized as follows:

\begin{itemize}
\item First, in the section \ref{sec:intro1} we describe the model, proposed
in \cite{ShelRuz2013}, and modify this model.

\item In section \ref{OneDim} we explain the main idea of the modification
in the Shelukhin-R\r{u}\v{z}i\v{c}ka model and collect some technical
results that are used in our main result, related with the proof of
global-in-time solvability result for the modified model. 
In particular, we introduce a potential for the viscous part of the stress
tensor and characterizes completely its sub-differential (see Proposition \ref%
{convexe1}). 

\item In section \ref{form} we formulate boundary-value problem and a
global-in-time existence theorem \ref{theorem 1} for the modified model.

\item In section \ref{Existence} we introduce an approximated problem %
\eqref{y1}, depending on a regularization index $n\in \mathbb{N}$, and show
the solvability of this approximated problem \eqref{y1} by Schauder fixed
point theorem. Also we derive a priori estimates for the solution of %
\eqref{y1}, which are independent on $n$;

\item Finally, in section \ref{Limit} we prove Theorem \ref{theorem 1},
applying the Lions-Aubin compactness theorem and a priori estimates of the
section \ref{Existence}.
\end{itemize}

\bigskip

\bigskip

\section{Model of asymmetric Bingham fluids}

\label{sec:intro1}

In what follows we explain the mathematical model of asymmetric Bingham
fluid that was proposed in the article \cite{ShelRuz2013}. For simplicity of
consideration in this article we consider a particular case when the angular
velocity field is zero.

For any matrix $X\in \mathbb{R}^{3\times 3}$ \ we define the symmetric and
asymmetric parts%
\begin{equation}
X_{s}=\frac{1}{2}\left( X+X^{T}\right) \qquad \text{and}\qquad X_{a}=\frac{1%
}{2}\left( X-X^{T}\right)  \label{s_a}
\end{equation}%
with the adjoint matrix $X^{T}$ having the property $\left( X^{T}\right)
_{ij}=x_{ji}$.\ \ Also we denote 
\begin{equation*}
X^{d}=X-\frac{\mbox{tr}\,X}{3}I\qquad \text{with}\qquad \mbox{tr}%
\,X=\sum_{i=1}^{3}x_{ii}.
\end{equation*}%
For any matrices $X,Y\in \mathbb{R}^{3\times 3}$ \ \ the scalar product $X:Y$
and the modulus $|X|$ of $X$ are defined by 
\begin{equation}
X:Y=\sum_{i,j=1}^{3}x_{ij}y_{ij},\qquad |X|=\left( X:X\right) ^{1/2}.
\label{dot}
\end{equation}

Let $\mathbf{v}=\mathbf{v}(\mathbf{x},t)$ be the velocity of the mass center
of the material point $({\boldsymbol{\xi }},t)$ for an asymmetric Bingham
fluid. We denote the rate of strain tensor 
\begin{equation}
\mathrm{B}=\mathrm{B}(\mathbf{v})=\frac{\partial \mathbf{v}}{\partial 
\mathbf{x}},\qquad (\partial \mathbf{v}/\partial \mathbf{x})_{ij}=\partial
v_{i}/\partial x_{j},  \label{ab}
\end{equation}%
and introduce the matrix 
\begin{equation*}
\mathrm{B}_{0}=\mathrm{B}_{s}+\varkappa \mathrm{B}_{a},\qquad \varkappa =%
\frac{2\mu _{1}}{\mu _{2}}
\end{equation*}%
where 
\begin{equation}
\qquad \mathrm{B}_{s}=\left( \frac{\partial \mathbf{v}}{\partial \mathbf{x}}%
\right) _{s}\qquad \text{and}\qquad \mathrm{B}_{a}=\left( \frac{\partial 
\mathbf{v}}{\partial \mathbf{x}}\right) _{a}  \label{ab00}
\end{equation}%
are the symmetric and asymmetric parts of $\mathrm{B}=\mathrm{B}(\mathbf{v}%
), $ respectively. \ The positive constants $\mu _{i}$ are viscosities of
the asymmetric Bingham fluid.\ \ An instant stress state of the fluid is
described by the Cauchy stress tensor $\mathrm{T}=-p\,\mathrm{I}+\mathrm{S}$%
, where $p$ and $\mathrm{S}$ are the pressure and the viscous part of the
stress tensor. In [19] the viscous part $\mathrm{S}$ of the stress tensor of
the fluid was suggested to be expressed as%
\begin{equation}
\mathrm{S}=\left\{ 
\begin{array}{ll}
2\mu _{1}\mathrm{B}_{0}+\tau _{\ast }\frac{\mathrm{B}_{0}}{|\mathrm{B}_{0}|},
& \mathrm{B}_{0}\neq 0 \\ 
\mathrm{S}_{plug}, & \mathrm{B}_{0}=0%
\end{array}%
\right.  \label{plug}
\end{equation}%
for some tensor $\ \mathrm{S}_{plug}\in \mathbb{R}^{3\times 3},$ such that$%
\quad |\mathrm{S}_{plug}|\leqslant \tau _{\ast }.$ The positive constant $%
\tau _{\ast }$ is the yield stress.

Finally, we have the momentum balance law describing the motion of
asymmetric fluid of Bingham 
\begin{equation}
\rho \dot{\mathbf{v}}=\mbox{div}\,\mathrm{T}+\rho \mathbf{f},
\label{momentum1}
\end{equation}%
where $\rho $ is the density and $\mathbf{f}$ is the mass force vector.

\bigskip

%

\begin{remark}
\label{rem1} Let us introduce the potential 
\begin{equation*}
V(X)=\mu _{1}|X|^{2}+\tau _{\ast }|X|,\qquad \forall \,X\in \mathbb{R}%
^{3\times 3}.
\end{equation*}%
Then the constitutive law \eqref{plug} for the asymmetric Bingham fluid can
be formulated as $\mathrm{S}\in \partial V(\mathrm{B}_{0}).$ \ Let us remind
that this inclusion is equivalent to the variational inequality 
\begin{equation}
V(X)-V(\mathrm{B}_{0})\geqslant \mathrm{S}:(X-\mathrm{B}_{0}),\qquad \forall
\,X\in \mathbb{R}^{3\times 3}.  \label{mon}
\end{equation}
\end{remark}

\bigskip

Nevertheless that the relation \eqref{plug} describes a plug zone in the
Bingham fluid, there exists a significant restriction in such modelling. As
we mentioned in \ref{Introduction} Introduction, in the articles \cite{smr}, 
\cite{ShelChem} the solvability of the model \eqref{plug}-\eqref{momentum1}
was shown only in the case when the asymmetric part is not present in the
tensor $\mathrm{S}$.\ As it is well known, one of the principal approach for
the study of problems with inclusions is the theory of monotone operators,
that was developed by Duvaut, Lions \cite{DL}.\ \ The major difficulty in
the study of Shelukhin-R\r{u}\v{z}i\v{c}ka model \cite{ShelRuz2013} consists
in the presence of the term $\mathrm{S}:\mathrm{B}_{0}$ in the inequality %
\eqref{mon}. The asymmetric part $\mathrm{B}_{a}$ of $\mathrm{B}_{0}$\ does
not permit to apply the theory of monotone operators.

In the following we present our modification in the above described model %
\eqref{plug}-\eqref{momentum1} in such way that permits to apply the theory
of monotone operators. Moreover, in our model the viscous part $\mathrm{S}$
\ \ will be a more general then in the mathematical model of \cite%
{ShelRuz2013}. For vectors functions $\mathbf{v}\in \mathbb{R}^{3}$ \ we
consider $\mathrm{B}=\mathrm{B}(\mathbf{v})$ introduced in \eqref{ab}. For
such defined tensor $\mathrm{B}$ \ we introduce the following tensors 
\begin{equation*}
\mathrm{B}_{\mu }=2\mu _{1}|\mathrm{B}_{s}|^{p-2}\mathrm{B}_{s}+\mu _{2}|%
\mathrm{B}_{a}|^{p-2}\mathrm{B}_{a},\qquad \mathrm{B}_{\nu }=|\mathrm{B}%
_{s}|^{\frac{p-2}{2}}\mathrm{B}_{s}+\nu |\mathrm{B}_{a}|^{\frac{p-2}{2}}%
\mathrm{B}_{a},\qquad 
\end{equation*}%
\begin{equation}
\mathrm{B}_{\nu ^{2}}=|\mathrm{B}_{s}|^{p-2}\mathrm{B}_{s}+\nu ^{2}|\mathrm{B%
}_{a}|^{p-2}\mathrm{B}_{a},  \label{ab0}
\end{equation}%
where $p\geqslant 2$, $\mu _{1},\mu _{2}$ are the viscosities of the
viscoplastic fluid of Bingham and $\tau _{\ast }$ \ is a so-called \textit{%
plug parameter}.\ Let us define the viscous part $\mathrm{S}$ of asymmetric
fluid of Bingham by 
\begin{equation}
\mathrm{S}=\left\{ 
\begin{array}{ll}
\mathrm{B}_{\mu }+\widehat{\tau }_{\ast }\frac{\mathrm{B}_{\nu ^{2}}}{|%
\mathrm{B}_{\nu }|^{\frac{2(p-1)}{p}}}, & \mathrm{B}_{\nu }\neq 0, \\ 
\mathrm{S}_{plug}, & \mathrm{B}_{\nu }=0,%
\end{array}%
\right.   \label{polarBingham1}
\end{equation}%
for some tensor $\ \mathrm{S}_{plug}=\mathrm{S}_{plug}(\mathbf{x},t)\in 
\mathbb{R}^{3\times 3},$ which fulfils the restriction $|\mathrm{S}%
_{plug}|\leqslant \tau _{\ast }$\textit{. }Here we denote 
\begin{equation}
\widehat{\tau }_{\ast }=\frac{\tau _{\ast }}{\max (1,\nu ^{2/p})}.
\label{tau*}
\end{equation}%
The major explanation of this modification is based on Proposition \ref%
{convexe1} proved in the following section.

\section{Some useful results}

\bigskip \label{OneDim}

In what follows the following algebraic result will be very useful.

\begin{lemma}
\label{lem1}\ The space of the matrices endowed with the dot product %
\eqref{dot} is the direct sum of the spaces of symmetric matrices and
anti-symmetric matrices. More precisely, for any matrices $X,$ $Y\in \mathbb{%
R}^{3\times 3},$\ we have%
\begin{align*}
& X_{s}:Y=X_{s}:Y_{s},\quad X_{a}:Y=X_{a}:Y_{a},\quad
X_{s}:Y_{a}=X_{a}:Y_{s}=0, \\
& X_{s}:X=X_{s}:X_{s}=\left\vert X_{s}\right\vert ^{2},\qquad
X_{a}:X=\left\vert X_{a}\right\vert ^{2},
\end{align*}%
with the notation introduced in \eqref{s_a}.
\end{lemma}

\textbf{Proof}. \ \ The coefficients of the matrices $X_{s}$ and $X_{a}$ are
equal to $x_{ij}^{s}=\frac{x_{ij}+x_{ji}}{2}$ and $x_{ij}^{a}=\frac{%
x_{ij}-x_{ji}}{2}$. Then \ $x_{ij}^{s}=x_{ji}^{s}$, $x_{ij}^{a}=-x_{ji}^{a},$%
\begin{equation*}
\begin{array}{ccl}
X_{s}:Y & = & \sum_{i,j=1}^{3}x_{ij}^{s}y_{ij}=\frac{1}{2}\left[
\sum_{i,j=1}^{3}x_{ij}^{s}y_{ij}+\sum_{i,j=1}^{3}x_{ji}^{s}y_{ji}\right] \\ 
& = & \sum_{i,j=1}^{3}x_{ij}^{s}\frac{1}{2}\left( y_{ij}+y_{ji}\right)
=\sum_{i,j=1}^{3}x_{ij}^{s}y_{ij}^{s}=X_{s}:Y_{s}%
\end{array}%
\end{equation*}%
and 
\begin{equation*}
\begin{array}{ccl}
X_{a}:Y & = & \sum_{i,j=1}^{3}x_{ij}^{a}y_{ij}=\frac{1}{2}\left[
\sum_{i,j=1}^{3}x_{ij}^{a}y_{ij}+\sum_{i,j=1}^{3}x_{ji}^{a}y_{ji}\right] \\ 
& = & \frac{1}{2}\left[ \sum_{i,j=1}^{3}x_{ij}^{s}y_{ij}-%
\sum_{i,j=1}^{3}x_{ij}^{a}y_{ji}\right] =\sum_{i,j=1}^{3}x_{ij}^{s}\frac{1}{2%
}\left( y_{ij}-y_{ji}\right) \\ 
& = & \sum_{i,j=1}^{3}x_{ij}^{a}y_{ij}^{a}=X_{a}:Y_{a}.%
\end{array}%
\end{equation*}%
Moreover, we have%
\begin{equation*}
\begin{array}{rcl}
\left( \frac{X\pm X^{T}}{2}\right) :\left( \frac{Y\mp Y^{T}}{2}\right) & = & 
\frac{1}{4}\sum_{i,j=1}^{3}(x_{ij}\pm x_{ji})(y_{ij}\mp y_{ji}) \\ 
& = & \frac{1}{4}\sum_{i,j=1}^{3}(x_{ij}y_{ij}-x_{ji}y_{ji})=0.%
\end{array}%
\end{equation*}%
$\hfill \;\blacksquare $

\bigskip

Before proceeding let us recall two basic theorems on convex analysis.

\begin{theorem}
\label{dir der} (see \cite[Theorem 23.1]{R} or \cite[Proposition 17.2]{BC})\
Let $f$ be a convex function from $\mathbb{R}^{n}$ to $[-\infty ,+\infty ]$,
and let $\mathbf{x}$ be a point where $f$ is finite. Then, for each $\mathbf{%
y}\in \mathbb{R}^{n}$ there exists the one-sided directional derivative of $%
f $ at $\mathbf{x}$ with respect to the vector $\mathbf{y}$, i.e. 
\begin{equation*}
f^{\prime }(\mathbf{x};\mathbf{y})=\lim_{\lambda \searrow 0}\frac{f(\mathbf{x%
}+\lambda \mathbf{y})-f(\mathbf{x})}{\lambda }.
\end{equation*}%
In fact, the difference quotient $(f(\mathbf{x}+\lambda \mathbf{y})-f(%
\mathbf{x}))/\lambda $ is a non-decreasing function of $\lambda >0$, so that 
\begin{equation*}
f^{\prime }(\mathbf{x};\mathbf{y})=\inf_{\lambda >0}(f(\mathbf{x}+\lambda 
\mathbf{y})-f(\mathbf{x}))/\lambda .
\end{equation*}
\end{theorem}

\begin{remark}
\label{hom} For any function $f:\mathbb{R}^{n}\rightarrow \mathbb{R}^{m}$
positively homogeneous of order 1, we have 
\begin{equation*}
f^{\prime }(0;\mathbf{y})=f(\mathbf{y}),\quad \forall \mathbf{y}\in \mathbb{R%
}^{n}.
\end{equation*}

Indeed, in this case one has 
\begin{equation*}
f^{\prime }(0;\mathbf{y})=\lim_{\lambda \searrow 0}\frac{f(\lambda \mathbf{y}%
)-f(0)}{\lambda }=\lim_{\lambda \downarrow 0}\frac{\lambda f(\mathbf{y})}{%
\lambda }=f(\mathbf{y}).
\end{equation*}
\end{remark}

Let us remember the concept of sub-differential.

\begin{definition}
A vector $\mathbf{x}^{\ast }$ is said to be a sub-gradient of $f$ at $%
\mathbf{x}$ if 
\begin{equation*}
f(\mathbf{y})\geqslant f(\mathbf{x})+(\mathbf{x}^{\ast },\mathbf{y}),\quad
\forall \mathbf{y}\in \mathbb{R}^{n}.
\end{equation*}%
The set of all sub-gradients of $f$ at $\mathbf{x}$ is called sub-differential
of $f$ \ at $\mathbf{x}$ and is denoted by $\partial f(\mathbf{x}).$
\end{definition}

\begin{theorem}
\label{dir} (see \cite[Theorem 23.2]{R} or \cite[Proposition 17.7]{BC}) Let $%
f$ be a convex function, and let $\mathbf{x}$ be a point where $f$ is
finite. Then $\mathbf{x}^{\ast }$ is a sub-gradient of $f$ at $\mathbf{x}$
if and only if 
\begin{equation*}
f^{\prime }(\mathbf{x};\mathbf{y})\geqslant (\mathbf{x}^{\ast },\mathbf{y}%
),\quad \forall \mathbf{y}\in \mathbb{R}^{n}.
\end{equation*}
\end{theorem}

\smallskip

Let us also remark the following fact.

\begin{remark}
\label{W} Let 
\begin{equation*}
\Vert x\Vert _{l^{p}(\mathbb{R}^{n})}=\sqrt[p]{|x_{1}|^{p}+...+|x_{n}|^{p}}%
,\qquad \mathbf{x}=(x_{1},\cdots ,x_{n})\in \mathbb{R}^{n},
\end{equation*}%
denote the $l^{p}$ norm in $\mathbb{R}^{n}$. Then the ${l}^{p}$ norm is
decreasing with respect to $p\in \lbrack 1,\infty ]$. \ 

This is easy to prove: given $p_{1}\leqslant p_{2}$ in $[1,\infty ]$ and $%
\mathbf{x}=(x_{1},\cdots ,x_{n})\not=0,$ let$\quad y_{i}=|x_{i}|/\Vert
x\Vert _{p_{2}}.$\ Then $|y_{i}|\leqslant 1$, so 
\begin{equation*}
|y_{i}|^{p_{1}}\geqslant |y_{i}|^{p_{2}},\qquad \Vert \mathbf{y}\Vert _{{l}%
^{p_{1}}}\geqslant 1,
\end{equation*}
and, consequently, $\Vert \mathbf{x}\Vert _{{l}^{p_{1}}}\geqslant \Vert 
\mathbf{x}\Vert _{{l}^{p_{2}}}$.
\end{remark}

\smallskip

\bigskip

\bigskip

The following results are used in a crucial way in the proof of our Theorem %
\ref{theorem 1}.

\begin{proposition}
\label{convexe1} For given $p\geqslant 2$ let us introduce the potential 
\begin{equation}
V(X)=\frac{2\mu _{1}}{p}|X_{s}|^{p}+\frac{\mu _{2}}{p}|X_{a}|^{p}+\widehat{%
\tau }_{\ast }||\mathrm{X}_{s}|^{\frac{p-2}{2}}\mathrm{X}_{s}+\nu |\mathrm{X}%
_{a}|^{\frac{p-2}{2}}\mathrm{X}_{a}|^{\frac{2}{p}}  \label{f1}
\end{equation}%
for any matrix $X\in \mathbb{R}^{3\times 3}.$ The potential $V$ is convex,
differentiable at any $X\not=0$ (equivalently, $X_{\nu }\not=0$) in $\mathbb{%
R}^{3\times 3}$, with 
\begin{equation*}
DV(X)=X_{\mu }+\widehat{\tau }_{\ast }\frac{X_{\nu ^{2}}}{|X_{\nu }|^{\frac{%
2(p-1)}{p}}},
\end{equation*}%
the matrices $X_{\mu },\mathrm{\ }X_{\nu },\mathrm{\ }X_{\nu ^{2}}$ are
defined in \eqref{ab0} (instead of $B$ we substitute $X$). \ \ Moreover:

$(a)\ \qquad B_{r_{p}}(0)\subset (\widehat{\tau }_{\ast })^{-1}\partial
V(0)\subset B_{\max \{1,\nu ^{\frac{2}{p}}\}}(0),$\quad where $B_{r}(0)$ is
the closed ball of a radius $r$ at the center $0$ and 
\begin{equation*}
r_{p}=\nu ^{2}/(1+\left( \nu ^{2}\right) ^{^{\frac{2}{p-2}}})^{^{\frac{p-2}{2%
}}};
\end{equation*}

$(b)\qquad $let \ \ $q=p/(p-1),$ then 
\begin{equation*}
(\widehat{\tau }_{\ast })^{-1}\partial V(0)=\{S\in \mathbb{R}^{3\times
3}:\quad |S_{s}|^{q}+\nu ^{2(1-q)}|S_{a}|^{q}\leqslant 1\}.
\end{equation*}
\end{proposition}

\textbf{Proof}. \ \ We write $V(X)=U(X)+\widehat{\tau }_{\ast }W(X)$ with%
\begin{equation*}
U(X)=\frac{2\mu _{1}}{p}|X_{s}|^{p}+\frac{\mu _{2}}{p}|X_{a}|^{p},\qquad
W(X)=||X_{s}|^{\frac{p-2}{2}}X_{s}+\nu |X_{a}|^{\frac{p-2}{2}}X_{a}|^{\frac{2%
}{p}}.
\end{equation*}

To see that $V$ is convex, we first notice that the functions 
\begin{equation*}
X\mapsto |X_{s}|^{p},\qquad X\mapsto |X_{a}|^{p}
\end{equation*}%
are convex since they are a composition of the convex function $t\in \mathbb{%
R}\mapsto |t|^{\frac{p}{2}}$ with quadratic functions $X\mapsto |X_{s}|^{2},$
$\ ~X\mapsto |X_{a}|^{2}$. \ In addition, we can check that the function $W$
is convex using that any norm is convex and the fact that 
\begin{equation*}
W(X)=\sqrt[p]{|X_{s}|^{p}+\nu ^{2}|X_{a}|^{p}}=\Vert (|X_{s}|,\nu ^{\frac{2}{%
p}}|X_{a}|)\Vert _{{l}^{p}(\mathbb{R}^{2})}
\end{equation*}%
is also a norm. Thus, $V$ is convex because it is a linear combination of
convex functions.

The function $V$ is differentiable at any $X\in \mathbb{R}^{3\times
3}\backslash \{0\}$ due to the chain rule, and we can compute $DV(X)$
differentiating directly $U$ and $W$ with respect to the standard variables $%
x_{ij}$ in $\mathbb{R}^{3\times 3}$, or, alternatively, using the chain
rule, one has%
\begin{align}
DV(X)=& 2\mu _{1}|X_{s}|_{s}^{p-2}X+\mu _{2}|X_{a}|_{a}^{p-2}X+\widehat{\tau 
}_{\ast }DW(X)  \notag \\
=& X_{\mu }+\widehat{\tau }_{\ast }(|X_{s}|^{p}+\nu ^{2}|X_{a}|^{p})^{\frac{1%
}{p}-1}(|X_{s}|^{p-2}X_{s}+\nu ^{2}|X_{a}|^{p-2}X_{a})  \notag \\
=& X_{\mu }+\widehat{\tau }_{\ast }|X_{\nu }|^{-\frac{2(p-1)}{p}}X_{\nu
}^{2}.  \label{DV}
\end{align}

Now, the function $U$ is differentiable also at $X=0$ and $DU(0)=0$. Then,
the sub-differential $\partial V(0)$ is equal to $\widehat{\tau }_{\ast
}\partial W(0)$. Let us show the items (a) and (b) in the statement of the Proposition).

By Remark \ref{hom}, 
\begin{equation*}
W^{\prime }(0;Y)=W(Y)\qquad \text{for any\quad }Y\in \mathbb{R}^{3\times 3}.
\end{equation*}%
Then, by Theorem \ref{dir}, 
\begin{equation}
S\in \partial W(0)=\widehat{\tau }_{\ast }^{-1}\partial V(0)\qquad
\Longleftrightarrow \qquad W(Y)\geqslant S:Y,\qquad \forall \,Y\in \mathbb{R}%
^{3\times 3}.  \label{r}
\end{equation}%
Taking in the inequality $Y=S$ and using the Remark \ref{W}, we obtain%
\begin{align*}
|S|^{2}& \leqslant \Vert (|S_{s}|,\nu ^{2/p}|S_{a}|)\Vert _{{l}^{p}(\mathbb{R%
}^{2})}\leqslant \max \{1,\nu ^{2/p}\}\Vert (|S_{s}|,|S_{a}|)\Vert _{{l}^{p}(%
\mathbb{R}^{2})} \\
& \leqslant \max \{1,\nu ^{2/p}\}\Vert (|S_{s}|,|S_{a}|)\Vert _{{l}^{2}(%
\mathbb{R}^{2})}=\max \{1,\nu ^{2/p}\}|S|,
\end{align*}%
hence we have proved the claim $\partial W(\mathrm{0})\subset %
B_{\max \{1,\nu ^{2/p}\}}(0)$ of item (a) of the Proposition.

Now let us show the claim 
\begin{equation}
B_{r_p}(0) \subset \partial W(\mathrm{0})=\widehat{%
\tau }_{\ast }^{-1}\partial V(0)  \label{ss}
\end{equation}%
of item (a). Let us consider arbitrary matrix $S\in \mathbb{R}^{3\times 3}$
satisfying the property 
\begin{equation}
|S||Y|\leqslant W(Y),\qquad \forall Y\in \mathbb{R}^{3\times 3}.  \label{s}
\end{equation}%
Then, by the Cauchy-Schwarz inequality, we have also 
\begin{equation*}
S:Y\leqslant W(Y),\qquad \forall \,Y\in \mathbb{R}^{3\times 3}.
\end{equation*}%
By Theorem \ref{dir} we conclude that any $S\in \mathbb{R}^{3\times 3},$
satisfying the property \eqref{s}, belongs to $\partial W(0)$. On the other
hand, by the positive homogeneity of $W$ the property \eqref{s} is
equivalent to 
\begin{equation*}
|S|\leqslant \min_{Y}\{W(Y):~Y\in \mathbb{R}^{3\times 3}\quad 
\mbox{   with 
}|Y|=1\}.
\end{equation*}%
Let us demonstrate that this minimum is equal to $r_{p}$, that gives the
claim \eqref{ss}. the exact value of $r_{p}$ is defined in the statement of
item (a) in the Proposition. Since $1=|Y|^{2}=|Y_{s}|^{2}+|Y_{a}|^{2}$, writing $%
t=|Y|_{a}^{2},$ we have $|Y_{s}|^{2}=1-t$ and%
\begin{equation*}
W(Y)\equiv W(t)=\sqrt[p]{(1-t)^{\frac{p}{2}}+\nu ^{2}t^{\frac{p}{2}}},\qquad
t\in \lbrack 0,1].
\end{equation*}%
By a straightforward computation, we obtain that the minimum of the function 
\begin{equation*}
\alpha (t)=(1-t)^{\frac{p}{2}}+\nu ^{2}t^{\frac{p}{2}}\qquad \text{in the
interval \quad }[0,1]
\end{equation*}%
is $\left( r_{p}\right) ^{p}$, and it is attained at $t_{\ast }=1/(1+\nu ^{%
\frac{4}{p-2}})$. Thus, we have proved the claim (a).

To show claim (b) we follow a similar reasoning as above, and with the help
of the H\"{o}lder inequality. Let $S\in \partial W(0)$. Then, by Theorem \ref%
{dir}, we have 
\begin{equation*}
S:Y\leqslant W(Y),\qquad \forall \,Y\in \mathbb{R}^{3\times 3}.
\end{equation*}%
If we take in this inequality $Y,$ having 
\begin{equation*}
Y_{s}=|S_{s}|^{q-2}S_{s}\qquad \text{and}\qquad Y_{a}=\nu ^{-\frac{2q}{p}%
}|S_{a}|^{q-2}S_{a},
\end{equation*}%
we obtain 
\begin{equation*}
|S_{s}|^{q}+\nu ^{-\frac{2q}{p}}|S_{a}|^{q}\leqslant \sqrt[p]{%
|S_{s}|^{q}+\nu ^{-2q+2}|S_{a}|^{q}}.
\end{equation*}%
This implies that 
\begin{equation*}
|S_{s}|^{q}+\nu ^{2(1-q)}|S_{a}|^{q}\leqslant 1,
\end{equation*}%
because $-2q/p=2(1-q)$ and $p\geqslant 2$. Reciprocally, if 
\begin{equation*}
|S_{s}|^{q}+\nu ^{2(1-q)}|S_{a}|^{q}\leqslant 1
\end{equation*}%
then, for all $Y\in \mathbb{R}^{3\times 3}$, using the H\"{o}lder
inequality, we have 
\begin{eqnarray*}
S &:&Y=S_{s}:Y_{s}+(\nu ^{-2/p}S_{a}):(\nu ^{2/p}Y_{a}) \\
&\leqslant &\sqrt[q]{|S_{s}|^{q}+\nu ^{-2q/p}|S_{s}|^{q}}\sqrt[p]{%
|Y_{s}|^{p}+\nu ^{2}|Y_{a}|^{p}}\leqslant W(Y).
\end{eqnarray*}%
Then, again by Theorem \ref{dir}, we obtain that $S\in \partial W(0)$. $%
\hfill \;\blacksquare $

\medskip 

\bigskip 

Now we show the auxiliary result that explains the definition of $\widehat{%
\tau }_{\ast }$ by \eqref{tau*} in \eqref{polarBingham1}$.$ 

\begin{corollary}
\label{est DW} For any matrix $B\in \mathbb{R}^{3\times 3}\backslash \{%
\mathrm{0}\}$ one has the estimate 
\begin{equation}
\frac{|B_{\nu }^{2}|}{\sqrt[p]{|B_{\nu }|^{2(p-1)}}}\leqslant \max \{1,\nu ^{%
\frac{2}{p}}\}.  \label{1c}
\end{equation}
\end{corollary}

\bigskip

\textbf{Proof}. \ \ As the derivative of the functional $W(X)$ has been
calculated in \eqref{DV} and equals to%
\begin{equation*}
DW(B)=\frac{B_{\nu }^{2}}{\sqrt[p]{|B_{\nu }|^{2(p-1)}}}\qquad \text{at any
\quad }B\not=0\text{.}
\end{equation*}

Now we claim that $DW(B)\in \partial W(0)$, and thus the estimate \eqref{1c}
shall follow from this fact, by (a) of Proposition \ref{convexe1}, since 
\begin{equation*}
\partial W(\mathrm{0})=(\widehat{\tau }_{\ast })^{-1}\partial V(\mathrm{0}%
)\subset B_{\max \{1,\nu ^{\frac{2}{p}}\}}(0).
\end{equation*}%
Accounting \eqref{r} we have to show the claim 
\begin{equation*}
DW(B):Y\leqslant W(Y)\qquad \text{for all}\quad Y\in \mathbb{R}^{3\times 3}.
\end{equation*}%
By Theorem \ref{dir der} we have that the directional derivative 
\begin{equation*}
DW(B):Y=\inf_{\lambda >0}\lambda ^{-1}\left( W(B+~\lambda Y)-W(B)\right) .
\end{equation*}%
In addition, $W$ is a norm by Remark \ref{W}, then 
\begin{equation*}
W(B+\lambda Y)\leqslant W(B)+\lambda W(Y).
\end{equation*}%
Summing up, we obtain our above claim.$\hfill \;\blacksquare $

\begin{remark}
\label{rem2} By Proposition \ref{convexe1}, the relation %
\eqref{polarBingham1} is equivalent to the variational inequality 
\begin{equation*}
V(X)-V(\mathrm{B})\geqslant S:(X-\mathrm{B}),\qquad \forall \,X\in \mathbb{R}%
^{3\times 3}.
\end{equation*}
\end{remark}

\bigskip

Next we present two technical results we shall use latter.

\begin{lemma}
Let $W=W(X)$ be a positive convex function on $X\in \mathbb{R}^{3\times 3}.$
\ Then, for any natural $n$, the approximated function 
\begin{equation*}
W_{n}(X)=\sqrt[p]{\left( W(X)\right) ^{p}+n^{-1}}
\end{equation*}%
is also convex with respect of the parameter $X\in \mathbb{R}^{3\times 3}$.
\end{lemma}

\bigskip \textbf{Proof}. \ \ Note that the function $\varphi (z)=\sqrt[p]{%
z^{p}+n^{-1}}$ is monotone increasing and convex function for $z\geqslant 0.$
\ Therefore \ applying the definition of convex function, we easily derive
that the composition $W_{n}(X)=\varphi (W(X))$ is also convex with respect
of the parameter $X\in \mathbb{R}^{3\times 3}$. $\hfill \;\blacksquare $

\bigskip

\begin{lemma}
\label{apr copy(1)} Let $n$ be an arbitrary natural number. We consider the
convex potential 
\begin{equation*}
V^{n}(X)=\frac{2\mu _{1}}{p}|X_{s}|^{p}+\frac{\mu _{2}}{p}|X_{a}|^{p}+%
\widehat{\tau }_{\ast }\sqrt[p]{||\mathrm{X}_{s}|^{\frac{p-2}{2}}\mathrm{X}%
_{s}+\nu |\mathrm{X}_{a}|^{\frac{p-2}{2}}\mathrm{X}_{a}|^{2}+n^{-1}},
\end{equation*}%
defined for arbitrary $X\in \mathbb{R}^{3\times 3}$. Let 
\begin{equation*}
\mathrm{S}^{n}=\mathrm{B}_{\mu}+\widehat{\tau }_{\ast }\frac{\mathrm{B}_{\nu
^{2}}}{\sqrt[p]{(|\mathrm{B}_{\nu }|^{2}+n^{-1})^{p-1}}}.
\end{equation*}
Then, for any\ given $B\in \mathbb{R}^{3\times 3}$, we have $\frac{\partial
V^{n}}{\partial X}(\mathrm{B})=\mathrm{S}^{n}$, i.e. 
\begin{equation*}
V^{n}(X)-V^{n}(\mathrm{B})\geqslant \mathrm{S}^{n}:(X-\mathrm{B}),\qquad
\forall \,X\in \mathbb{R}^{3\times 3}.
\end{equation*}
\end{lemma}

\bigskip \textbf{Proof}. \ \ Straightforward computation. $\hfill
\;\blacksquare $

\section{Statement of the problem}

\bigskip

\label{form}

\bigskip

Let us consider the motion of an asymmetric Bingham fluid, assuming that the
fluid is incompressible. For simplicity of considerations we admit that the
density $\rho $ is equal to $1$ and 
neglect the mass force vector$\ \mathbf{f}$. \ Then \ the flow equations %
\eqref{momentum1} for the velocity $\mathbf{v}$\ in a bounded domain $\Omega
\subset \mathbb{R}^{3}$\ with the boundary $\Gamma $\ are 
\begin{equation}
\mathbf{v}_{t}+\left( \mathbf{v\cdot \nabla }\right) \mathbf{v}=\mbox{div}\ 
\mathrm{T},\qquad \mbox{div}\,\mathbf{v}=0\qquad \text{in \ \ }\Omega
_{T}=(0,T)\times \Omega ,  \label{OneDeq1}
\end{equation}%
where $\mathrm{T}=-p\,\mathrm{I}+\mathrm{S}$ and $\mathrm{S}$ satisfies the
constitutive law \eqref{polarBingham1} with the relations \eqref{ab}, %
\eqref{ab00}, \eqref{ab0}.\ \ \ We add to this system the initial data%
\begin{equation}
\mathbf{v}|_{t=0}=\mathbf{v}_{0}\qquad \text{in}\quad \Omega .
\label{dataOneDim1}
\end{equation}

The system \eqref{OneDeq1} is mostly supplemented with the usual Dirichlet
boundary condition. The Dirichlet condition implies the adherence of fluid
particles to the boundary. For the motion of Bingham fluids (such as the
extrusion of the toothpaste from the tube, the mud flows in drilling
engineering, \ the propane flow within the hydro-fracture, etc.) \ it is
more natural to study the system \eqref{OneDeq1} with slip type boundary
conditions, permitting the slippage of the fluid against the boundary.\ To
describe accurately physical phenomena we consider homogeneous \ Navier slip
boundary conditions%
\begin{equation}
\mathbf{v}\cdot \mathbf{n}=0,\;\quad \left[ \mathrm{T}\,\mathbf{n}+\alpha 
\mathbf{v}\right] \cdot {\bm{\tau }}=0\qquad \text{on }\ \Gamma
_{T}=(0,T)\times \Gamma ,  \label{Direchlet}
\end{equation}%
where $\alpha $ is a positive \textit{friction} coefficient. For the
discussion of the Navier slip boundary conditions we refer to the articles 
\cite{CC1}-\cite{CC4}.

\bigskip

Let us introduce some notations to formulate our main result. \ We denote by 
$(\cdot ,\cdot )$ the 
dot product in $L^{2}(\Omega )$. Also we define the spaces%
\begin{eqnarray}
H &=&\{\mathbf{v}\in L^{2}(\Omega ):\,\mbox{div }\mathbf{v}=0\;\ \text{ in}\;%
\mathcal{D}^{\prime }(\Omega ),\;\ \mathbf{v}\cdot \mathbf{n}=0\;\ \text{ in}%
\;H^{-1/2}(\Gamma )\},  \notag \\
V &=&\{\mathbf{v}\in H^{1}(\Omega ):\,\mbox{div }\mathbf{v}=0\;\ \text{ a.e.
in \ }\Omega ,\;\ \mathbf{v}\cdot \mathbf{n}=0\;\ \text{ in}\;H^{1/2}(\Gamma
)\},  \notag \\
V_{p} &=&\{\mathbf{v}\in V:\,\;\;|\nabla \mathbf{v}|^{p}\in L^{1}(\Omega )\}.
\label{w}
\end{eqnarray}%
The space $V_{p}$ is endowed 
with the norm $\Vert \mathbf{v}\Vert _{V_{p}}=\Vert \mathbf{v}\Vert
_{L^{2}(\Omega )}+\Vert \nabla \mathbf{v}\Vert _{L^{p}(\Omega )}.$

The main objective of our article is to show the well-posedness of the
system \eqref{OneDeq1}-\eqref{Direchlet} for unknown functions $\mathbf{v}$
and $\mathrm{S}$. This result of the well-posedness is formulated in the
following theorem, in which we also define the concept of the weak solution
for the system \eqref{OneDeq1}-\eqref{Direchlet}. This concept is a direct
consequence of the equations \eqref{OneDeq1}, the boundary conditions %
\eqref{Direchlet} and the integral equality 
\begin{equation}
-\int_{\Omega }\mbox{div}\,\mathrm{T}\cdot {\boldsymbol{\varphi }}\ d\mathbf{%
x}=-\int_{\Gamma }\left( \mathrm{T}\ \mathbf{n}\right) \cdot {\boldsymbol{%
\varphi }}\ d\mathbf{\gamma }+\int_{\Omega }\mathrm{T}:\frac{\partial {%
\boldsymbol{\varphi }}}{\partial \mathbf{x}}\,d\mathbf{x},  \label{yi3}
\end{equation}%
which is valid for any $(3\times 3)$-matrix function $\mathrm{T}\in
H^{1}(\Omega )$\ and any 3D-vector function ${\boldsymbol{\varphi }}\in
H^{1}(\Omega ).$

\begin{theorem}
\label{theorem 1} 
Let $\Omega$ be a bounded domain in $\mathbb{R}^3$ with a $C^{1}-$smooth
boundary $\Gamma$, $p\geqslant 2$, a given real number, 
\begin{equation}
\mathbf{v}_{0}\in H\qquad \mbox{ and } \qquad \alpha \in L^{2}(0,T;L^{\infty
}(\Omega )).  \label{in}
\end{equation}%
Then there exists a function $\mathbf{v}$ and a $(3\times 3)-$matrix
function $\mathrm{S}$, such that 
\begin{equation}
\mathbf{v}\in L^{\infty }(0,T;H)\cap L^{2}(0,T;V_{p}),\quad \quad \mathbf{v}%
_{t}\in L^{2}(0,T;V_{p}^{\ast })  \label{regvw}
\end{equation}%
and the pair $\left( \mathbf{v,\ }\mathrm{S}\right) $ is a weak solution of
the system \eqref{OneDeq1}-\eqref{Direchlet}, satisfying\ the integral
equality 
\begin{equation}
\int\limits_{\Omega _{T}}\left[ \mathbf{v}\partial _{t}{\boldsymbol{\varphi }%
}+\left( \mathbf{v\otimes v}-\mathrm{S}\right) :\frac{\partial {\boldsymbol{%
\varphi }}}{\partial \mathbf{x}}\right] \,d\mathbf{x}dt+\int\limits_{\Omega }%
\mathbf{v}_{0}{\boldsymbol{\varphi }}(0)\,d\mathbf{x}=\int_{\Gamma
_{T}}\alpha (\mathbf{v}\cdot {\boldsymbol{\varphi }})\,d\mathbf{\gamma }\,dt
\label{weakForm}
\end{equation}%
for any test function 
\begin{equation}
{\boldsymbol{\varphi }}\in H^{1}(0,T;V_{p}) \quad \mbox{such that } \quad {%
\boldsymbol{\varphi }} (\cdot ,T)=0.  \label{a}
\end{equation}%
The $(3\times 3)-$matrix function $\ \mathrm{S}\in L^{p/(p-1)}(\Omega _{T})$
fulfils the relation \eqref{polarBingham1}.

Moreover, if 
\begin{equation*}
p\geqslant \frac{7+\sqrt{19}}{5}\approx 2.272,
\end{equation*}
then the solution $\left( \mathbf{v,\ }\mathrm{S}\right) $ is unique.
\end{theorem}

\bigskip

\bigskip

\section{Construction of the approximation problem}

\label{Existence}

\bigskip

In this section we consider an approximated problem for the system %
\eqref{OneDeq1}-\eqref{Direchlet} and 
solve this approximated problem applying the Faedo-Galerkin method and the
Schauder fixed point argument (see for instance \cite{AKM}).\ 

\bigskip

Since the space $V_{p}$ is separable, it is the span of a countable set of
linearly independent functions $\left\{ \mathbf{e}_{k}\right\}
_{k=1}^{\infty }$. 
More precisely, we can choose this set as the eigenfunctions for the
following non-linear Stokes type equations with Navier boundary conditions: 
\begin{equation*}
\left\{ 
\begin{array}{ll}
-{\mbox{div}}(\mathrm{T}(\mathbf{e}_{k}))=\lambda _{k}\mathbf{e}_{k},\qquad 
\mathrm{div}\,\mathbf{e}_{k}=0 & \quad \mbox{in}\ \Omega ,\vspace{2mm} \\ 
&  \\ 
\mathbf{e}_{k}\cdot \mathbf{n}=0,\;\quad \left[ \mathrm{T}(e_{k})\,\mathbf{n}%
+\alpha \mathbf{e}_{k}\right] \cdot {\bm{\tau }}=0 & \quad \mbox{on}\ \Gamma%
\end{array}%
\right.
\end{equation*}%
with \ \ $\mathrm{T}(\mathbf{e}_{k})=-p_{k}\,\mathrm{I}+|\nabla \mathbf{e}%
_{k}|^{p-2}\nabla \mathbf{e}_{k}.$ The solvability of this problem follows
from the spectral theory \cite{Evans}.\ \ This theory permits to construct
this set $\left\{ \mathbf{e}_{k}\right\} _{k=1}^{\infty }$ as an orthogonal
basis for $V_{p}$ and an orthonormal basis for $H.$

We can consider the subspace $V_{p}^{n}=\mathrm{span}\,\{\mathbf{e}%
_{1},\ldots ,\mathbf{e}_{n}\}$ of $V_{p}$, \ 
for any fixed natural $n.$ \ Let us define the vector function%
\begin{equation}
\mathbf{v}^{n}(t)=\sum_{k=1}^{n}c_{k}^{(n)}(t)\ \mathbf{e}_{k},\qquad
c_{k}^{(n)}(t)\in {\mathbb{R}},  \label{vv}
\end{equation}%
as the solution of the approximate system%
\begin{equation}
\left\{ 
\begin{array}{l}
\int\limits_{\Omega }[\partial _{t}\mathbf{v}^{n}\mathbf{e}_{k}+\left( 
\mathbf{v}^{n}\mathbf{\cdot \nabla }\right) \mathbf{v}^{n}\mathbf{e}_{k}+%
\mathrm{S}^{n}:\frac{\partial \mathbf{e}_{k}}{\partial \mathbf{x}}]\,d%
\mathbf{x}+\int_{\Gamma }\alpha (\mathbf{v}^{n}\cdot {\bm{\tau }})(\mathbf{e}%
_{k}\cdot {\bm{\tau }})\,d\mathbf{\gamma }\,dt=0, \\ 
\qquad \forall k=1,2,\dots ,n, \\ 
\mathbf{v}^{n}(0)=\mathbf{v}_{0}^{n}.%
\end{array}%
\right.  \label{y1}
\end{equation}%
The $(3\times 3)-$matrix functions $\mathrm{T}^{n},$ $\ \ \mathrm{S}^{n}$\
are prescribed by the relations%
\begin{equation}
\mathrm{T}^{n}=-p^{n}\,\mathrm{I}+\mathrm{S}^{n},\qquad \mathrm{S}^{n}=%
\mathrm{B}_{\mu }^{n}+\widehat{\tau }_{\ast }\frac{\mathrm{B}_{\nu ^{2}}^{n}%
}{\sqrt[p]{(|\mathrm{B}_{\nu }^{n}|^{2}+n^{-1})^{p-1}}}  \label{delta}
\end{equation}%
and the matrix functions $\mathrm{B}^{n}=\mathrm{B}(\mathbf{v}^{n}),$ $%
\mathrm{B}_{\mu }^{n},\mathrm{\ B}_{\nu }^{n},\mathrm{\ B}_{\nu ^{2}}^{n}$
are calculated through the formulas \eqref{ab}, \eqref{ab00}, \eqref{ab0}. \
The function $\mathbf{v}_{0}^{n}$ \ is the orthogonal projection of $\mathbf{%
v}_{0}\in H$\ into the space $V_{p}^{n}.$ Note that the system \eqref{vv}-%
\eqref{delta} is a weak formulation of the problem%
\begin{equation*}
\left\{ 
\begin{array}{c}
\mathbf{v}_{t}^{n}+\left( \mathbf{v}^{n}\mathbf{\cdot \nabla }\right) 
\mathbf{v}^{n}=\mbox{div}\ \mathrm{T}^{n},\qquad \mbox{div}\,\mathbf{v}%
^{n}=0,\qquad \text{in}\quad \Omega _{T}, \\ 
\mathbf{v}^{n}\cdot \mathbf{n}=0,\;\quad \left[ \mathrm{T}^{n}\mathbf{n}%
+\alpha \mathbf{v}^{n}\right] \cdot {\bm{\tau }}=0\qquad \text{on }\ \Gamma
_{T}, \\ 
\mathbf{v}^{n}|_{t=0}=\mathbf{v}_{0}^{n}\qquad \text{in}\quad \Omega .%
\end{array}%
\right.
\end{equation*}%
Next we prove that the approximated problem \eqref{y1} is solvable.

\begin{lemma}
\label{apr} Let us assume that the data $\mathbf{v}_{0},$ $\alpha $ \
satisfy the conditions \eqref{in}. Then there exists a solution $\mathbf{v}%
^{n}\in L^{\infty }(0,T;V_{p})$ \ \ of the system \eqref{vv}-\eqref{delta},
satisfying the a priori estimate%
\begin{equation}
\int\limits_{\Omega }|\mathbf{v}^{n}|^{2}d\mathbf{x}+\int\limits_{0}^{t}%
\left[ \int\limits_{\Omega }|\frac{\partial \mathbf{v}^{n}}{\partial \mathbf{%
x}}|^{p}d\mathbf{x}+\int_{\Gamma }\alpha |\mathbf{v}^{n}|^{2}\,d\mathbf{%
\gamma }\right] dt\leqslant A,\quad t\in \lbrack 0,T],  \label{EST}
\end{equation}%
%
%
%
%
%
%
%
%
%
%
%
%
%
%
%
%
%
%
%
%
%
%
%
%
%
%
where $A$ is a constant independent on $n$. More precisely, $A$ depends only
on the data $\mathbf{v}_{0},\mu _{i},\nu$.
\end{lemma}

\textbf{Proof}.\ \ The system \eqref{y1} is a system of $n$ ordinary
differential equations of the first order, which can be written in the form%
\begin{equation*}
\frac{d\mathbf{c}^{(n)}}{dt}=F(\mathbf{c}^{(n)}),\qquad t\in (0,T),
\end{equation*}%
for the vector function $\mathbf{c}%
^{(n)}(t)=(c_{1}^{(n)}(t),...,c_{n}^{(n)}(t))$ with $c_{k}^{(n)}$ introduced
in \eqref{vv}. 
We can solve system using the Schauder fixed point theorem.

For an arbitrary $\widehat{\mathbf{c}}^{(n)}\in C([0,T])$, \ we define 
\begin{equation*}
\widehat{\mathbf{v}}^{n}(t)=\sum_{k=1}^{n}\widehat{c}_{k}^{(n)}(t)\ \mathbf{e%
}_{k}\qquad \text{and}\quad \widehat{\tau }_{\ast }(\widehat{\mathbf{c}}%
^{(n)})=\frac{\mathrm{1}}{\sqrt[p]{(|\mathrm{B}_{\nu }^{n}(\widehat{\mathbf{v%
}}^{n})|^{2}+n^{-1})^{p-1}}}.
\end{equation*}%
Let $\mathbf{c}^{(n)}$ be an unknown, such that the vector function%
\begin{equation}
\mathbf{v}^{n}(t)=\sum_{k=1}^{n}c_{k}^{(n)}(t)\ \mathbf{e}_{k}  \label{bs0}
\end{equation}%
solves the system%
\begin{equation}
\left\{ 
\begin{array}{l}
\int\limits_{\Omega }[\partial _{t}\mathbf{v}^{n}\mathbf{e}_{k}+\left( 
\widehat{\mathbf{v}}^{n}\mathbf{\cdot \nabla }\right) \mathbf{v}^{n}\mathbf{e%
}_{k}+\widehat{\mathrm{S}}^{n}:\frac{\partial \mathbf{e}_{k}}{\partial 
\mathbf{x}}]\,d\mathbf{x}+\int_{\Gamma }\alpha (\mathbf{v}^{n}\cdot {%
\bm{\tau }})(\mathbf{e}_{k}\cdot {\bm{\tau }})\,d\mathbf{\gamma }\,dt=0, \\ 
\qquad \forall k=1,2,\dots ,n,\vspace{2mm} \\ 
\mathbf{v}^{n}(0)=\mathbf{v}_{0}^{n}%
\end{array}%
\right.  \label{bs3}
\end{equation}%
with the $(3\times 3)-$matrix function%
\begin{equation}
\widehat{\mathrm{S}}^{n}=\mathrm{B}_{\mu }^{n}+\widehat{\tau }_{\ast }(%
\widehat{\mathbf{c}}^{(n)})\mathrm{B}_{\nu ^{2}}^{n},  \label{bs}
\end{equation}%
and the matrix functions $\mathrm{B}^{n}=\mathrm{B}(\mathbf{v}^{n}),$ $%
\mathrm{B}_{\mu }^{n},\mathrm{\ B}_{\nu }^{n},\mathrm{\ B}_{\nu ^{2}}^{n}$ 
given by the formulas \eqref{ab}, \eqref{ab00} and \eqref{ab0}. From the
theory of ordinary differential equations, it follows that the linear system %
\eqref{bs0}-\eqref{bs3}, of $n$ ordinary linear differential equations, has
an unique solution $\mathbf{c}^{(n)}\in C^{1}([0,T])$. Therefore, we can
consider the operator $K:C([0,T])\rightarrow C([0,T])$ defined as 
\begin{equation*}
\mathbf{c}^{(n)}=K\left( \widehat{\mathbf{c}}^{(n)}\right) .
\end{equation*}%
The solvability of the system \eqref{vv}-\eqref{delta} will be shown if we
demonstrate that this operator $K$ has a fixed point, which we shall do by
the Schauder fixed point theorem. Thus, we have to prove that this operator
is compact on a bounded convex subset $M$ \ of $\ C([0,T])$.

First, let us deduce a priori estimates for $\mathbf{c}^{(n)}.$ \ We
multiply \eqref{bs3}$_{1}$ by $c_{k}^{(n)}$ and take the sum on the index $%
k=1,...,n.$\ Then the integration over the time interval $(0,t)$ gives%
\begin{equation}
\frac{1}{2}\int_{\Omega }|\mathbf{v}^{n}|^{2}d\mathbf{x}+\int_{0}^{t}\left[
\int_{\Omega }\widehat{\mathrm{S}}^{n}:\frac{\partial \mathbf{v}^{n}}{%
\partial \mathbf{x}}d\mathbf{x}+\int_{\Gamma }\alpha |\mathbf{v}^{n}|^{2}\,d%
\mathbf{\gamma }\right] \,dt=\frac{1}{2}\int_{\Omega }|\mathbf{v}%
_{0}^{n}|^{2}d\mathbf{x}.  \label{c}
\end{equation}%
Lemma \ref{lem1} and the definition \eqref{bs} of $\widehat{\mathrm{S}}^{n}$
\ imply

\begin{equation*}
\mathrm{B}_{\mu }^{n}:\frac{\partial \mathbf{v}^{n}}{\partial \mathbf{x}}%
=2\mu _{1}|\mathrm{B}_{s}^{n}|^{p}+\mu _{2}|\mathrm{B}_{a}^{n}|^{p},\qquad 
\mathrm{B}_{\nu ^{2}}^{n}:\frac{\partial \mathbf{v}^{n}}{\partial \mathbf{x}}%
=|\mathrm{B}_{s}^{n}|^{p}+\nu ^{2}|\mathrm{B}_{a}^{n}|^{p}
\end{equation*}%
and 
\begin{equation*}
\widehat{\mathrm{S}}^{n}:\frac{\partial \mathbf{v}^{n}}{\partial \mathbf{x}}%
=2\mu _{1}|\mathrm{B}_{s}^{n}|^{p}+\mu _{2}|\mathrm{B}_{a}^{n}|^{p}+\widehat{%
\tau }_{\ast }(\widehat{\mathbf{c}}^{(n)})\left( |\mathrm{B}%
_{s}^{n}|^{p}+\nu ^{2}|\mathrm{B}_{a}^{n}|^{p}\right) .
\end{equation*}%
Therefore, we deduce the inequality%
\begin{eqnarray}
\frac{1}{2}\int_{\Omega }|\mathbf{v}^{n}|^{2}d\mathbf{x}&+&\int_{0}^{t}\left[
\int_{\Omega }\left( 2\mu _{1}|\mathrm{B}_{s}^{n}|^{p}+\mu _{2}|\mathrm{B}%
_{a}^{n}|^{p}\right) d\mathbf{x}+\int_{\Gamma }\alpha |\mathbf{v}^{n}|^{2}\,d%
\mathbf{\gamma }\right]  \notag \\
&\leqslant& \frac{1}{2}\int_{\Omega }|\mathbf{v}_{0}^{n}|^{2}d\mathbf{x,}
\label{r0}
\end{eqnarray}%
which gives the apriori estimate \eqref{EST}.

Also, we have 
\begin{equation*}
|\mathrm{B}_{\mu }^{n}|^{2}=4\mu _{1}^{2}|\mathrm{B}_{s}^{n}|^{2(p-1)}+\mu
_{2}^{2}|\mathrm{B}_{a}^{n}|^{2(p-1)},\qquad |\mathrm{B}_{\nu
^{2}}^{n}|^{2}=|\mathrm{B}_{s}^{n}|^{2(p-1)}+\nu ^{4}|\mathrm{B}%
_{a}^{n}|^{2(p-1)},
\end{equation*}%
by Lemma \ref{lem1}. Since 
\begin{equation*}
\widehat{\tau }_{\ast }(\widehat{\mathbf{c}}^{(n)})\leqslant \frac{\mathrm{1}%
}{\sqrt[p]{n^{-(p-1)}}},\quad \text{then \ \ }|\widehat{\mathrm{S}}%
^{n}|\leqslant |\mathrm{B}_{\mu }^{n}|+\frac{\mathrm{1}}{\sqrt[p]{n^{-(p-1)}}%
}|\mathrm{B}_{\nu ^{2}}^{n}|,
\end{equation*}%
and using \eqref{r0}, we obtain%
\begin{equation}
||\widehat{\mathrm{S}}^{n}||_{L^{p/(p-1)}(\Omega _{T})}\leqslant C(n)
\label{EST1}
\end{equation}%
where the constant $C(n)$\ depends only on $n$.

Let us define the bounded convex set 
\begin{equation*}
M=\left\{ \mathbf{c}^{(n)}\in C([0,T]):||\mathbf{c}^{(n)}||\leqslant
A\right\} ,
\end{equation*}%
where the norm $||\mathbf{c}^{(n)}||^{2}=\max_{t\in \lbrack
0,T]}\sum\limits_{k=1}^{n}\left( c_{k}^{(n)}(t)\right) ^{2}$ is defined on
the space $C([0,T])$ and the constant $A$ is prescribed in \eqref{EST}. \ 

Let us assume that $\widehat{\mathbf{c}}^{(n)}\in M$. \ Since $\{\mathbf{e}%
_{j}\}_{j=1}^{\infty }$ is the orthonormal basis\ for the space $H,$ \ then
the equality \eqref{bs3} can be written as%
\begin{equation*}
\frac{dc_{k}^{(n)}}{dt}=-\int\limits_{\Omega }\left[ \left( \widehat{\mathbf{%
v}}^{n}\mathbf{\cdot \nabla }\right) \mathbf{v}^{n}\mathbf{e}_{k}+\widehat{%
\mathrm{S}}^{n}:\frac{\partial \mathbf{e}_{k}}{\partial \mathbf{x}}\right]
\,d\mathbf{x}-\int_{\Gamma }\alpha (\mathbf{v}^{n}\cdot {\bm{\tau }})(%
\mathbf{e}_{k}\cdot {\bm{\tau }})\,d\mathbf{\gamma }\,dt
\end{equation*}%
for any $k=1,\dots ,n.$ \ Since $V_{p}\subset H^{1}(\Omega ),$ then the
Sobolev continuous embedding$\ H^{1}(\Omega )\hookrightarrow L^{6}(\Omega
)\cap L^{2}(\Gamma )$ and the Holder inequality imply%
\begin{eqnarray}
|\frac{dc_{k}^{(n)}}{dt}| &\leqslant &\{\Vert \widehat{\mathbf{v}}^{n}\Vert
_{L^{4}(\Omega )}\Vert \nabla \mathbf{v}^{n}\Vert _{L^{2}(\Omega )}+\Vert 
\sqrt{\alpha }\Vert _{L^{\infty }(\Omega )}\Vert \sqrt{\alpha }\mathbf{v}%
^{n}\Vert _{L^{2}(\Gamma )}\}\Vert \mathbf{e}_{k}\Vert _{H^{1}(\Omega )} 
\notag \\
&&\quad +||\widehat{\mathrm{S}}^{n}||_{L^{p/(p-1)}(\Omega )}\Vert \frac{%
\partial \mathbf{e}_{k}}{\partial \mathbf{x}}\Vert _{L^{p}(\Omega )},
\label{ccc}
\end{eqnarray}%
that is%
\begin{eqnarray}
|\frac{dc_{k}^{(n)}}{dt}| &\leqslant & C(n)\{\Vert \widehat{\mathbf{v}}%
^{n}\Vert _{L^{4}(\Omega )}\Vert \nabla \mathbf{v}^{n}\Vert _{L^{2}(\Omega
)}+\Vert \sqrt{\alpha }\mathbf{v}^{n}\Vert _{L^{2}(\Gamma )}  \notag \\
&& \qquad + ||\widehat{\mathrm{S}}^{n}||_{L^{p/(p-1)}(\Omega )}\}.
\label{c_k}
\end{eqnarray}%
Let us recall the Gagliardo--Nirenberg-Sobolev inequality (see \cite{tem})%
\begin{align}
&\int_{0}^{T}\left( \Vert \widehat{\mathbf{v}}^{n}\Vert _{L^{4}(\Omega
)}\Vert \nabla \mathbf{v}^{n}\Vert _{L^{2}(\Omega )}\right) ^{8/7}\,dt
\leqslant C\int_{0}^{T}\left( ||\widehat{\mathbf{v}}^{n}||_{L^{2}(\Omega
)}^{1/4}||\nabla \mathbf{v}^{n}||_{L^{2}(\Omega )}^{7/4}\right) ^{8/7}\,dt 
\notag \\
&\qquad \qquad\qquad \leqslant C\Vert \widehat{\mathbf{v}}^{n}\Vert
_{L^{\infty }\left( 0,T;L^{2}(\Omega )\right) }^{2/7}\Vert \mathbf{v}%
^{n}\Vert _{L^{2}\left( 0,T;H^{1}(\Omega )\right) }^{2}.  \label{vn}
\end{align}%
Therefore, the integration of the inequality \eqref{c_k} over the time
interval $[t,t+h]\subset \lbrack 0,T]$, the H\"older inequality and $%
\widehat{\mathbf{c}}^{(n)}\in M$, plus the estimates \eqref{EST}, \eqref{vn}%
, imply 
\begin{equation*}
|c_{k}^{(n)}(t+h)-c_{k}^{(n)}(t)|\leqslant \int_{t}^{t+h}|\frac{dc_{k}^{(n)}%
}{dt}|\leqslant C(n)h^{q}\qquad \text{with\quad }q=min(\frac{1}{8},\frac{1}{2%
},\frac{1}{p})
\end{equation*}%
for each $k=1,\dots ,n.$ Therefore the operator $K:M\rightarrow M$ is
compact by the Arzela-Ascoli theorem.

The continuity of $K$ is a direct consequence of the theorem on continuous
dependence of the solution of the Cauchy problem \eqref{bs0}-\eqref{bs3}
with respect to the coefficients $\widehat{\mathbf{c}}^{(n)}$ . \ Therefore,
the operator $K$\ fulfils the conditions of the Schauder fixed point
theorem, which implies the existence of a fixed point of $K,$ and gives the
solution of the system \eqref{vv}-\eqref{delta}. $\hfill \;\blacksquare $

\bigskip

\begin{lemma}
\label{1} We assume that \ the data $\mathbf{v}_{0},$ $\alpha $ \ fulfil
the conditions \eqref{in}. Then there exists a solution 
\begin{equation*}
\mathbf{v}^{n}\in L^{\infty }(0,T;H)\cap L^{2}(0,T;V_{p})
\end{equation*}%
of the system \eqref{vv}-\eqref{delta} that satisfies the following
estimates 
\begin{equation}
\int\limits_{\Omega }|\mathbf{v}^{n}|^{2}d\mathbf{x}+\int\limits_{0}^{t}%
\left[ \int\limits_{\Omega }|\frac{\partial \mathbf{v}^{n}}{\partial \mathbf{%
x}}|^{p}d\mathbf{x}+\int_{\Gamma }\alpha |\mathbf{v}^{n}|^{2}\,d\mathbf{%
\gamma }\right] dt\leqslant C,\quad t\in \lbrack 0,T],  \label{omegadelta}
\end{equation}%
\begin{equation}
||\mathrm{S}^{n}||_{L^{p/(p-1)}(\Omega _{T})}\leqslant C  \label{Sdelta}
\end{equation}%
and%
\begin{equation}
||\partial _{t}\mathbf{v}^{n}||_{L^{8/7}(0,T;\,\,V_{p}^{\ast })}\leqslant C.
\label{timedelta}
\end{equation}%
Here and below $C$ is a positive constant that \textit{does not depend on} $%
n,$ but may depend on $\mathbf{v}_{0}$ and $\alpha $.
\end{lemma}

\textbf{Proof}.\ \ Let us note that in Lemma \ref{apr} we already shown that
the solution $\mathbf{v}^{n}\in L^{\infty }(0,T;V_{p}^{n})$\ of \eqref{vv}-%
\eqref{delta}\ satisfies \eqref{omegadelta}. Moreover, from the definition %
\eqref{delta} of $\mathrm{S}^{n}$, we have%
\begin{equation*}
|\mathrm{S}^{n}|\leqslant |\mathrm{B}_{\mu }^{n}|+\tau _{\ast }.
\end{equation*}%
This is a direct consequence of the estimate \eqref{1c}.\ \ Hence we obtain
the a priori estimate \eqref{Sdelta} from \eqref{omegadelta} and the equality 
\begin{equation*}
|\mathrm{B}_{\mu }^{n}|^{2}=4\mu _{1}^{2}|\mathrm{B}_{s}^{n}|^{2(p-1)}+\mu
_{2}^{2}|\mathrm{B}_{a}^{n}|^{2(p-1)}.
\end{equation*}

\bigskip

Let us consider a subspace $V_{p}^{n}$\ of of $V_{p},$\ defined in Lemma \ref%
{apr}.\ Let $P_{n}$ be the orthogonal projection of $V_{p}$ \ onto $%
V_{p}^{n}.$\ Let ${\boldsymbol{\varphi }}\in H^{1}(0,T;V_{p})$ be an
arbitrary function.\ \ The first equality of \eqref{y1} is linear with
respect of the functions $\mathbf{e}_{k}$, $k=1,..n,$ then we have%
\begin{equation}
\left\{ 
\begin{array}{l}
\int\limits_{\Omega }[\partial _{t}\mathbf{v}^{n}(P_{n}{\boldsymbol{\varphi }%
})+\left( \mathbf{v}^{n}\mathbf{\cdot }\nabla \right) \mathbf{v}^{n}\ (P_{n}{%
\boldsymbol{\varphi }})+\mathrm{S}^{n}:\frac{\partial (P_{n}{\boldsymbol{%
\varphi }})}{\partial \mathbf{x}}]\,d\mathbf{x} \\ 
\qquad \qquad \qquad +\int_{\Gamma }\alpha (\mathbf{v}^{n}\cdot (P_{n}{%
\boldsymbol{\varphi }}))\,d\mathbf{\gamma }\,dt=0, \\ 
\mathbf{v}^{n}(0)=\mathbf{v}_{0}^{n},%
\end{array}%
\right.  \label{vnn}
\end{equation}%
since $\{\mathbf{e}_{j}\}_{j=1}^{\infty }$ \ is the orthogonal basis\ for
the space $V_{p}.$ \ 

As it was done in \eqref{ccc}, we obtain \ 
\begin{eqnarray*}
|(\partial _{t}\mathbf{v}^{n},{\boldsymbol{\varphi }})_{L_{2}(\Omega )}|
&\leqslant &\{\Vert \mathbf{v}^{n}\Vert _{L^{4}(\Omega )}\Vert \nabla 
\mathbf{v}^{n}\Vert _{L^{2}(\Omega )} \\
&+&\Vert \sqrt{\alpha }\Vert _{L^{\infty }(\Omega )}\Vert \sqrt{\alpha }%
\mathbf{v}^{n}\Vert _{L^{2}(\Gamma )}\}\Vert {P_{n}\boldsymbol{\varphi }}%
\Vert _{H^{1}(\Omega )} \\
&+&||\mathrm{S}^{n}||_{L^{p/(p-1)}(\Omega )}\Vert \frac{\partial ({P_{n}%
\boldsymbol{\varphi }})}{\partial \mathbf{x}}\Vert _{L^{p}(\Omega )},
\end{eqnarray*}%
then 
\begin{align}
& ||\partial _{t}\mathbf{v}^{n}||_{V_{p}^{\ast }}=\sup_{{\boldsymbol{\varphi 
}}\in V_{p}}\left\{ |(\partial _{t}\mathbf{v}^{n},{\boldsymbol{\varphi }}%
)_{L^{2}(\Omega )}|:\quad ||{\boldsymbol{\varphi }}||_{V_{p}}=1\right\} 
\notag \\
& \leqslant C\left\{ \Vert \mathbf{v}^{n}\Vert _{L^{4}(\Omega )}\Vert \nabla 
\mathbf{v}^{n}\Vert _{L^{2}(\Omega )}+\Vert \sqrt{\alpha }\mathbf{v}%
^{n}\Vert _{L^{2}(\Gamma )}+||\mathrm{S}^{n}||_{L^{p/(p-1)}(\Omega
)}\right\} ,  \label{as}
\end{align}%
since the norm $\Vert {\boldsymbol{\varphi }}\Vert _{V_{p}}=\Vert {%
\boldsymbol{\varphi }}\Vert _{L^{2}(\Omega )}+\Vert {\boldsymbol{\varphi }}%
\Vert _{L^{p}(\Omega )}$ on the space $V_{p}$ \ of the continuous operator $%
P_{n}$ is less or equal than $1$.

Therefore this last inequality \eqref{as}, the a priori estimates %
\eqref{omegadelta}-\eqref{Sdelta} and the Gagliardo--Nirenberg-Sobolev
inequality \eqref{vn}, written for $\mathbf{v}^{n}$ instead of \ $\widehat{%
\mathbf{v}}^{n}$ \ imply \ the estimate (\ref{timedelta}).$\hfill
\;\blacksquare $

\bigskip

\section{Limit transition}

\label{Limit}

We have that $V_{p}\subset H^{1}(\Omega ),$ then the embedding result $%
H^{1}(\Omega )\hookrightarrow L^{2}(\Gamma )$ and the estimates %
\eqref{omegadelta}-\eqref{timedelta} imply the existence of sub-sequences,
such that%
\begin{eqnarray}
\mathbf{v}^{n} &\rightharpoonup &\mathbf{v}\quad \mbox{ *-weakly in }\quad
L^{\infty }(0,T;L^{2}(\Omega )),  \notag \\
\mathbf{v}^{n} &\rightharpoonup &\mathbf{v}\quad \mbox{ weakly in }\quad
L^{2}(0,T;V_{p}\cap L^{2}(\Gamma )),  \notag \\
\mathrm{S}^{n} &\rightharpoonup &\mathrm{S}\quad \mbox{ weakly in}\quad
L^{2}(\Omega _{T}).  \label{conv}
\end{eqnarray}%
Let us remember the Aubin-Lions-Simon compactness result \cite{sim}, \cite%
{tem}.

\begin{lemma}
\label{ALS} Let $X_{0}$, $X$ and $X_{1}$ be three Banach spaces with $%
X_{0}\subseteq X\subseteq X_{1}$. Suppose that \bigskip $X_{0},$ $X_{1} $
are reflexive, $X_{0}$ is compactly embedded in $X$ and that $X$ is
continuously embedded in $X_{1}$. Let%
\begin{equation*}
W=\left\{ v\in L^{2}(0,T;X_{0}),\qquad \partial _{t}v\in
L^{8/7}(0,T;X_{1})\right\} .
\end{equation*}%
Then the embedding of $W$ into $L^{2}(0,T;X)$ is compact.
\end{lemma}

\bigskip

\bigskip

Therefore the estimates \eqref{omegadelta}, \eqref{timedelta}, the compact
embedding $H^{1}(\Omega )\hookrightarrow L^{2}(\Omega )$ and Lemma \ref{ALS}
give that%
\begin{equation}
\mathbf{v}^{n}\rightarrow \mathbf{v}\quad \mbox{strongly
in }~L^{2}(\Omega _{T})\quad \mbox{ and
a.e. in }~\Omega _{T}.  \label{convergence}
\end{equation}%
Hence \ applying the convergences \eqref{conv}-\eqref{convergence} in %
\eqref{vnn}, we deduce that the limit functions $\mathbf{v},$ $\mathrm{S}$\
fulfil the integral equality%
\begin{equation}
\int\limits_{\Omega _{T}}\left[ \mathbf{v}\partial _{t}{\boldsymbol{\varphi }%
}+\left( \mathbf{v\otimes v}-\mathrm{S}\right) :\frac{\partial {\boldsymbol{%
\varphi }}}{\partial \mathbf{x}}\right] \,d\mathbf{x}dt+\int\limits_{\Omega }%
\mathbf{v}_{0}{\boldsymbol{\varphi }}(0)\,d\mathbf{x}=\int_{\Gamma
_{T}}\alpha (\mathbf{v}\cdot {\boldsymbol{\varphi }})\,d\mathbf{\gamma }\,dt
\label{weak}
\end{equation}%
for any function ${\boldsymbol{\varphi }}\in H^{1}(0,T;V_{p}),\quad \,{%
\boldsymbol{\varphi }}(\cdot ,T)=0.$

In what follows we use the approach of the theory of variational
inequalities \cite{DL}, \cite{Evans} to demonstrate the relation %
\eqref{polarBingham1}, that ends the proof of Theorem \ref{theorem 1}. \ For
a fixed natural $n$\ we consider the convex potential 
\begin{equation*}
V^{n}=V^{n}(X),\forall \,X\in \mathbb{R}^{3\times 3},
\end{equation*}%
introduced in Lemma \ref{apr copy(1)}. By this Lemma we have 
\begin{equation}
\mathrm{S}^{n}=\frac{\partial V^{n}}{\partial X}(\mathrm{B}^{n})\qquad \text{%
with~ }\mathrm{B}^{n}=\frac{\partial \mathbf{v}^{n}}{\partial \mathbf{x}}
\label{s10}
\end{equation}%
and 
\begin{equation}
V^{n}(X)-V^{n}(\mathrm{B}^{n})\geqslant \mathrm{S}^{n}:(X-\mathrm{B}%
^{n})\quad \text{a.e. in }\Omega _{T},\quad \forall \,X\in \mathbb{R}%
^{3\times 3}.  \label{s0}
\end{equation}%
are equivalent. \ 

Let us denote $\Omega _{r}=(0,r)\times \Omega $ \ for arbitrary $r\in (0,T),$
the equality \eqref{c} can be written as%
\begin{equation}
\frac{1}{2}\int\limits_{\Omega }\left[ |\mathbf{v}^{n}|^{2}(r)-|\mathbf{v}%
_{0}^{n}|^{2}\right] \,d\mathbf{x}+\int_{0}^{r}\int_{\Gamma }\alpha |\mathbf{%
v}^{n}|^{2}\,d\mathbf{\gamma }\,dt=-\int\limits_{\Omega _{r}}\left[ \mathrm{S%
}^{n}:\frac{\partial \mathbf{v}^{n}}{\partial \mathbf{x}}\right] \,\,d%
\mathbf{x}\,dt.  \label{s2}
\end{equation}%
If we substitute \eqref{s2} in \eqref{s0}, then the lower semi-continuity
property of convex functional with respect of weak convergence gives%
\begin{eqnarray*}
&&\int\limits_{\Omega _{r}}V(X)-V(\mathrm{B})\,d\mathbf{x}dt\geqslant
\lim_{n\rightarrow \infty }\inf \int\limits_{\Omega _{r}}V^{n}(X)-V^{n}(%
\mathrm{B}^{n})\,d\mathbf{x}dt \\
&\geqslant &\lim_{n\rightarrow \infty }\inf \{\int\limits_{\Omega _{r}}%
\mathrm{S}^{n}:X\,d\mathbf{x}dt+\frac{1}{2}\int\limits_{\Omega }\left[ |%
\mathbf{v}^{n}|^{2}(r)-|\mathbf{v}_{0}^{n}|^{2}\right] \,d\mathbf{x}%
+\int_{0}^{r}\int_{\Gamma }\alpha |\mathbf{v}^{n}|^{2}\,d\mathbf{\gamma }%
\,dt\} \\
&\geqslant &\int\limits_{\Omega _{r}}\mathrm{S}:X\,d\mathbf{x}dt+\frac{1}{2}%
\int\limits_{\Omega }\left[ |\mathbf{v}|^{2}(r)-|\mathbf{v}_{0}|^{2}\right]
\,d\mathbf{x}+\int_{0}^{r}\int_{\Gamma }\alpha |\mathbf{v}|^{2}\,d\mathbf{%
\gamma }\,dt
\end{eqnarray*}%
by use the convergences \eqref{conv}. Hence for any matrix function $X\in
L^{2}(\Omega _{T})$ we have the inequality%
\begin{eqnarray}
\int\limits_{\Omega _{r}}V(X) &-&V(\mathrm{B})\,d\mathbf{x}dt\geqslant
\int\limits_{\Omega _{r}}\mathrm{S}:X\,d\mathbf{x}dt  \notag \\
&+&\frac{1}{2}\int\limits_{\Omega }\left[ |\mathbf{v}|^{2}(r)-|\mathbf{v}%
_{0}|^{2}\right] \,d\mathbf{x}+\int_{0}^{r}\int_{\Gamma }\alpha |\mathbf{v}%
|^{2}\,d\mathbf{\gamma }\,dt.  \label{eq0.0}
\end{eqnarray}

Let us take $\mathbf{\varphi }=\mathbf{v}(1-sgn_{+}^{\varepsilon }(t-r))$\
in the equality \eqref{weakForm} for a fixed $r\in (0,T),$ where%
\begin{equation*}
sgn_{+}^{\varepsilon }(t)=\left\{ 
\begin{array}{cc}
0, & \text{if \ }t<0; \\ 
t/\varepsilon , & \text{if \ }0\leqslant t<\varepsilon ; \\ 
1, & \text{if \ }\varepsilon \leqslant t.%
\end{array}%
\right.
\end{equation*}%
In the obtained equality the limit transition on $\varepsilon \rightarrow 0$
implies%
\begin{equation}
\frac{1}{2}\int\limits_{\Omega }\left[ |\mathbf{v}|^{2}(r)-|\mathbf{v}%
_{0}|^{2}\right] \,d\mathbf{x}+\int_{0}^{r}\int_{\Gamma }\alpha |\mathbf{v}%
|^{2}\,d\mathbf{\gamma }\,dt=-\int\limits_{\Omega _{r}}\left[ \mathrm{S}:%
\frac{\partial \mathbf{v}}{\partial \mathbf{x}}\right] \,\,d\mathbf{x}dt.
\label{eq00}
\end{equation}%
Substituting \eqref{eq00} in \eqref{eq0.0}, we derive%
\begin{equation*}
\int\limits_{\Omega _{r}}V(X)-V(\mathrm{B})\,d\mathbf{x}dt\geqslant
\int\limits_{\Omega _{r}}\mathrm{S}:(X-\mathrm{B})\,d\mathbf{x}dt\qquad 
\text{with \ }\mathrm{B}=\frac{\partial \mathbf{v}}{\partial \mathbf{x}}.
\end{equation*}%
Since the matrix function $X\in L^{2}(\Omega _{T})$ is arbitrary, we can
choose in this inequality $X=\mathrm{B}+\varepsilon \mathrm{Z}$ \ for any
positive $\varepsilon $ and any matrix function \textrm{$Z$}$\in
L^{2}(\Omega _{T}),$ \ which gives 
\begin{equation*}
\lim_{\varepsilon \rightarrow 0^{+}}\int\limits_{\Omega _{r}}\frac{V(\mathrm{%
B}+\varepsilon \mathrm{Z})-V(\mathrm{B})}{\varepsilon }\,d\mathbf{x}%
dt\geqslant \int\limits_{\Omega _{r}}\mathrm{S}:\mathrm{Z}\,d\mathbf{x}dt.
\end{equation*}%
So, passing the limit over the sign of integration, we obtain 
\begin{equation*}
\int\limits_{\Omega _{r}}V^{\prime }(\mathrm{B};\mathrm{Z})\,d\mathbf{x}%
dt\geqslant \int\limits_{\Omega _{r}}\mathrm{S}:\mathrm{Z}\,d\mathbf{x}dt,
\end{equation*}%
for any matrix function $\text{$Z$}\in L^{2}(\Omega _{T})$. Now, since $%
V^{\prime }(\mathrm{B};\mathrm{Z})$ and $\mathrm{S}:\mathrm{Z}$ are positively
homogeneous with respect to $Z$, this implies indeed in 
\begin{equation*}
\int\limits_{\Omega _{r}}V^{\prime }(\mathrm{B};\mathrm{Z})\xi \,d\mathbf{x}%
dt\geqslant \int\limits_{\Omega _{r}}(\mathrm{S}:\mathrm{Z})\xi \,d\mathbf{x}%
dt,
\end{equation*}%
for any positive function $\xi \in L^{\infty }(\Omega _{T})$ and any matrix
function $Z\in L^{2}(\Omega _{T})$. Therefore, $V^{\prime }(\mathrm{B};%
\mathrm{Z})\geqslant \mathrm{S}:\mathrm{Z}$, for any matrix function $%
\mathrm{Z}\in L^{2}(\Omega _{T})$, and thus, $S(x,t)\in \partial V(B(x,t))$
by Theorem \ref{dir} and $S$ has the form \eqref{polarBingham1} by
Proposition \ref{convexe1}. In particular, if $B=0$, we have $%
|S(x,t)|\leqslant \tau _{\ast }$ by (a) of Proposition \ref{convexe1}.

\bigskip

In the sequel we will show the uniqueness result.\ First let us show that
the norms 
\begin{equation*}
\Vert \mathbf{v}\Vert _{V_{p}}=\Vert \mathbf{v}\Vert _{L^{2}(\Omega )}+\Vert
\nabla \mathbf{v}\Vert _{L^{p}(\Omega )},\qquad \Vert \mathbf{v}\Vert
_{W_{p}^{1}}=\Vert \mathbf{v}\Vert _{L^{p}(\Omega )}+\Vert \nabla \mathbf{v}%
\Vert _{L^{p}(\Omega )}
\end{equation*}%
are equivalent. It is enough to show the following result.

\begin{lemma}
\label{equivalents} There exists a constant $C$, such that%
\begin{equation}
\Vert \mathbf{v}\Vert _{W_{p}^{1}(\Omega )}\leqslant C\Vert \mathbf{v}\Vert
_{V_{p}},\quad \forall \mathbf{v}\in W_{p}^{1}(\Omega ).  \label{eq}
\end{equation}
\end{lemma}

\textbf{Proof.} \ Assume that the affirmation of Lemma is not true, then for
any $n\in \mathbb{N}$ there exists a vector function $\mathbf{v}_{n}\in
W_{p}^{1}(\Omega ),$ such that $\Vert \mathbf{v}_{n}\Vert _{W_{p}^{1}(\Omega
)}>n\Vert \mathbf{v}_{n}\Vert _{V_{p}}.$ Let us define $\widetilde{\mathbf{v}%
}_{n}=\frac{\mathbf{v}_{n}}{\Vert \mathbf{v}\Vert _{W_{p}^{1}(\Omega )}},$
that fulfils%
\begin{equation}
\Vert \widetilde{\mathbf{v}}_{n}\Vert _{W_{p}^{1}(\Omega )}=1,\qquad \Vert 
\widetilde{\mathbf{v}}_{n}\Vert _{V_{p}}<\frac{1}{n}.  \label{mn}
\end{equation}%
Hence $\{\widetilde{\mathbf{v}}_{n}\}_{n=1}^{\infty }$ is compact in $%
L^{p}(\Omega ),$ such that there exists a strongly convergent subsequence $%
\widetilde{\mathbf{v}}_{n^{\prime }}$ to some $\widetilde{\mathbf{v}}\in
L^{p}(\Omega ).$ This subsequence $\widetilde{\mathbf{v}}_{n^{\prime }}$ is
a strongly convergent to $\widetilde{\mathbf{v}}$ in the space $L^{2}(\Omega
)$ by the inequality 
\begin{equation*}
\Vert \mathbf{v}\Vert _{L^{2}(\Omega )}\leqslant C\Vert \mathbf{v}\Vert
_{L^{p}(\Omega )},\quad \forall \mathbf{v}\in L^{p}(\Omega ).
\end{equation*}%
Also from \eqref{mn} we have that 
\begin{equation*}
\lim_{n\rightarrow \infty }\Vert \nabla \widetilde{\mathbf{v}}_{n}\Vert
_{L^{p}(\Omega )}=0,\quad \lim_{n\rightarrow \infty }\Vert \widetilde{%
\mathbf{v}}_{n}\Vert _{L^{p}(\Omega )}=1\qquad \text{and}\qquad
\lim_{n\rightarrow \infty }\Vert \widetilde{\mathbf{v}}_{n}\Vert
_{L^{2}(\Omega )}=0,
\end{equation*}%
that implies 
\begin{equation*}
\Vert \widetilde{\mathbf{v}}\Vert _{L^{p}(\Omega )}=1\qquad \text{and}\qquad
\Vert \widetilde{\mathbf{v}}\Vert _{L^{2}(\Omega )}=0,
\end{equation*}%
which is impossible. Hence \eqref{eq} is true.$\hfill \;\blacksquare $

\bigskip

Since for any $\beta \in [0,1]$, and any $a,b\geqslant 0$, we have $\left(
a+b\right)^{\beta }\leqslant a^{\beta }+b^{\beta }$, then Lemma \ref%
{equivalents} and the well known interpolation inequality 
$\Vert \mathbf{v}\Vert _{L^{r}(\Omega )}\leqslant C||\mathbf{v}%
||_{L^{2}(\Omega )}^{1-\beta }||\mathbf{v}||_{W_{p}^{1}(\Omega )}^{\beta }$, 
valid for any $\mathbf{v}\in W_{p}^{1}(\Omega )$ (see e.g. of \cite[Lemma 2.2%
]{AKM}), give the following interpolation result.

\begin{lemma}
\label{interpolation} There exists a positive constant $C$, such that%
\begin{equation}
\Vert \mathbf{v}\Vert _{L^{r}(\Omega )}\leqslant C(||\mathbf{v}%
||_{L^{2}(\Omega )}^{1-\beta }||\nabla \mathbf{v}||_{L^{p}(\Omega )}^{\beta
}+||\mathbf{v}||_{L^{2}(\Omega )}),\quad \forall \mathbf{v}\in V_{p},
\label{interp}
\end{equation}%
where $\beta =\left( \frac{1}{2}-\frac{1}{r}\right) /(\frac{5}{6}-\frac{1}{p}%
)$.
\end{lemma}



\bigskip

Now we are able to prove the uniqueness result. Let us denote the difference
of two functions $\mathbf{f}_{1}$ and $\mathbf{f}_{2}$ by $\overline{\mathbf{%
f}}$, i.e. $\overline{\mathbf{f}}=\mathbf{f}_{1}-\mathbf{f}_{2}$. \ Let us
admit the existence of two different solutions $\mathbf{v}_{1},\mathbf{v}%
_{2} $\ \ with respective tensors $\mathrm{S}_{1},$ $\mathrm{S}_{2},$\
satisfying the relation \eqref{polarBingham1}. \ By \eqref{weakForm} the
difference $\overline{\mathbf{v}}$ fulfils the equality%
\begin{equation}
\int\limits_{\Omega _{T}}\left[ \overline{\mathbf{v}}\partial _{t}{%
\boldsymbol{\varphi }}+\left( \overline{\mathbf{v\otimes v}}-\overline{%
\mathrm{S}}\right) :\frac{\partial {\boldsymbol{\varphi }}}{\partial \mathbf{%
x}}\right] \,d\mathbf{x}dt=\int_{\Gamma _{T}}\alpha (\overline{\mathbf{v}}%
\cdot {\boldsymbol{\varphi }})\,d\mathbf{\gamma }\,dt  \label{v1}
\end{equation}%
for any ${\boldsymbol{\varphi }}\in H^{1}(0,T;V_{p}),$ such that ${%
\boldsymbol{\varphi }}(\cdot ,T)=0$ in $\Omega .$

It is easily to check that%
\begin{equation*}
\int\limits_{\Omega }\left[ \left( \overline{\mathbf{v\otimes v}}\right) :%
\frac{\partial \overline{\mathbf{v}}}{\partial \mathbf{x}}\right] \,d\mathbf{%
x}=-\int\limits_{\Omega }\left[ \left( \overline{\mathbf{v}}\mathbf{\otimes }%
\overline{\mathbf{v}}\right) :\frac{\partial \mathbf{v}_{2}}{\partial 
\mathbf{x}}\right] \,d\mathbf{x.}\ 
\end{equation*}%
Also there exists a constant $\widetilde{C}>0,$ depending only on $\mu
_{1},\mu _{2}$ and $p,$ such that 
\begin{equation*}
C(|\overline{\mathrm{B}_{s}}|^{p}+|\overline{\mathrm{B}_{a}}|^{p})\leqslant 
 \overline{\mathrm{B}_{\mu }}:\frac{\partial \overline{%
\mathbf{v}}}{\partial \mathbf{x}}\leqslant \overline{\mathrm{S}}:\frac{%
\partial \overline{\mathbf{v}}}{\partial \mathbf{x}},
\end{equation*}%
being a consequence of the monotonicity of the second term in the relation %
\eqref{polarBingham1} and the inequality $(1.25)_{1}$ of Lemma 1.19, shown in 
\cite{MRN}.

Let us fix an arbitrary $r\in (0,T)$ and take $\mathbf{\varphi }=\overline{%
\mathbf{v}}(1-sgn_{+}^{\varepsilon }(t-r))$\ in \eqref{v1}. Then the limit
transition on $\varepsilon \rightarrow 0$\ in the obtained equality and the
H\"older inequality give that%
\begin{eqnarray}
&&\frac{1}{2}\int\limits_{\Omega }|\overline{\mathbf{v}}|^{2}(r)\,d\mathbf{x}%
+\widetilde{C}\int\limits_{\Omega _{r}}|\frac{\partial \overline{\mathbf{v}}%
}{\partial \mathbf{x}}|^{p}\,\,d\mathbf{x}dt  \notag \\
&\leqslant &\frac{1}{2}\int\limits_{\Omega }|\overline{\mathbf{v}}|^{2}(r)\,d%
\mathbf{x}+\int\limits_{\Omega _{r}}\left[ \overline{\mathrm{S}}:\frac{%
\partial \overline{\mathbf{v}}}{\partial \mathbf{x}}\right] \,\,d\mathbf{x}%
dt+\int_{\Gamma _{T}}\alpha |\overline{\mathbf{v}}|^{2}\,d\mathbf{\gamma }%
\,dt  \notag \\
&=&\int\limits_{\Omega _{r}}(\overline{\mathbf{v}}\mathbf{\otimes }\overline{%
\mathbf{v}}):\frac{\partial \mathbf{v}_{2}}{\partial \mathbf{x}}\,\,d\mathbf{%
x}dt.  \label{d1}
\end{eqnarray}%
In the sequel we follows the ideas presented in Theorem 3.2 of \cite{tem}
and Theorem 4.29 of \cite{MRN}. \ By H\"older's inequality and Lemma \ref%
{interpolation}, the right hand side of \eqref{d1} is estimated as 
\begin{eqnarray}
|\int\limits_{\Omega }(\overline{\mathbf{v}}\mathbf{\otimes }\overline{%
\mathbf{v}}) &:&\frac{\partial \mathbf{v}_{2}}{\partial \mathbf{x}}\,\,d%
\mathbf{x}|\leqslant C||\nabla \mathbf{v}_{2}||_{L^{p}(\Omega )}||\overline{%
\mathbf{v}}||_{L^{\frac{2p}{p-1}}(\Omega )}^{2}  \notag \\
&\leqslant &C||\nabla \mathbf{v}_{2}||_{L^{p}(\Omega )}\left( ||\overline{%
\mathbf{v}}||_{L^{2}(\Omega )}^{\frac{5p-9}{5p-6}}||\frac{\partial \overline{%
\mathbf{v}}}{\partial \mathbf{x}}||_{L^{p}(\Omega )}^{\frac{3}{5p-6}}+||%
\overline{\mathbf{v}}||_{L^{2}(\Omega )}\right) ^{2}  \notag \\
&\leqslant &C||\nabla \mathbf{v}_{2}||_{L^{p}(\Omega )}\left( ||\overline{%
\mathbf{v}}||_{L^{2}(\Omega )}^{2\left( \frac{5p-9}{5p-6}\right) }||\frac{%
\partial \overline{\mathbf{v}}}{\partial \mathbf{x}}||_{L^{p}(\Omega )}^{%
\frac{6}{5p-6}}+||\overline{\mathbf{v}}||_{L^{2}(\Omega )}^{2}\right)  \notag
\\
&\leqslant &\varepsilon ||\frac{\partial \overline{\mathbf{v}}}{\partial 
\mathbf{x}}||_{L^{p}(\Omega )}^{p}+C_{\varepsilon }||\nabla \mathbf{v}%
_{2}||_{L^{p}(\Omega )}^{\left( \frac{5p-6}{5p-9}\right) }||\overline{%
\mathbf{v}}||_{L^{2}(\Omega )}^{2}  \notag \\
&&+C||\nabla \mathbf{v}_{2}||_{L^{p}(\Omega )}||\overline{\mathbf{v}}%
||_{L^{2}(\Omega )}^{2},  \label{conv1}
\end{eqnarray}%
where the $\varepsilon -$version of Young's inequality has been used in the
last inequality. Hence taking $\varepsilon =\widetilde{C},$ we obtain%
\begin{equation*}
z(t)\leqslant C+C\int\limits_{0}^{t}f(s)z(s)\ ds\text{\qquad with\quad }%
z(t)=\int\limits_{\Omega }|\overline{\mathbf{v}}|^{2}\ d\mathbf{x,}\text{%
\quad }f(s)=||\nabla \mathbf{v}_{2}||_{L^{p}(\Omega )}^{\left( \frac{5p-6}{%
5p-9}\right) }.
\end{equation*}%
Therefore, if $p\geqslant \left( \frac{5p-6}{5p-9}\right) ,$\ that is $%
p\geqslant \frac{7+\sqrt{19}}{5}\approx 2.272,$ then applying the Gronwall
inequality we obtain $z(t)=0$ a.e. in $(0,T)$ and we deduce the
global-in-time uniqueness result.$\hfill \;\blacksquare $

\end{document}